\documentstyle[11pt]{article}

\font\twelvebf=cmbx12
\font\twelvebf=cmbx12

\newcommand{\ncm}{\newcommand}
\ncm{\rncm}{\renewcommand}
\rncm{\sec}{\setc{0}\section}
\ncm{\bsn}{\bigskip\noindent}
\ncm{\beq}{\begin{equation}}
\ncm{\eeq}{\end{equation}}
\ncm{\bea}{\begin{eqnarray}}
\ncm{\beanon}{\begin{eqnarray*}}
\ncm{\eea}{\end{eqnarray}}
\ncm{\eeanon}{\end{eqnarray*}}
\ncm{\fns}{\footnotesize}

\rncm{\theequation}{\thesection.\arabic{equation}}
\ncm{\setc}[1]{\setcounter{equation}{#1}}
\newcounter{eqnr}
\newcounter{axiom}

\newenvironment{axioms}{\stepcounter{axiom}
    \setcounter{eqnr}{\value{equation}}\setc{0}
    \rncm{\theequation}{A.\arabic{axiom}\alph{equation}}
    \begin{eqnarray}}{\end{eqnarray}\setc{\value{eqnr}}}

\newenvironment{axiom}{
   \setcounter{eqnr}{\value{equation}}
   \setc{\value{axiom}}
   \rncm{\theequation}{A.\arabic{equation}}
   \begin{equation}}
  {\end{equation}\setc{\value{eqnr}}}

\newenvironment{eqnarrayabc}{\stepcounter{equation}
  \setcounter{eqnr}{\value{equation}}\setc{0}
  \rncm{\theequation}{\thesection.\arabic{eqnr}\alph{equation}}
  \begin{eqnarray}}{\end{eqnarray}\setc{\value{eqnr}}}

\ncm{\bealph}{\begin{eqnarrayabc}}
\ncm{\eealph}{\end{eqnarrayabc}}
\ncm{\bit}{\begin{itemize}}
\ncm{\eit}{\end{itemize}}
\ncm{\eqboxabc}[3]{\newline\parbox[t]{1.5cm}{#1}\hfill
  \parbox[b]{12cm}{\begin{eqnarray*} #3\end{eqnarray*}}\hfill
   \parbox[b]{1.5cm}{\vspace{-0.0cm}
  \begin{eqnarrayabc}#2\end{eqnarrayabc}}\newline}

\newtheorem{thm}{Theorem}[section]
\newtheorem{prop}[thm]{Proposition}

\newtheorem{defi}[thm]{Definition}
\newtheorem{coro}[thm]{Corollary}

\def\lcros{\raise1.5pt\hbox{$\scriptstyle\triangleright$}\!
           \raise1.9pt\hbox{$\scriptscriptstyle < \,$}}
\def\truesupset{{\lower5pt\hbox{$\scriptstyle\supset$}\atop
 \raise5pt\hbox{$\scriptscriptstyle\not=$}}}
\def\truesubset{{\lower5pt\hbox{$\scriptstyle\subset$}\atop
 \raise5pt\hbox{$\scriptscriptstyle\not=$}}}
\def\End{\hbox{End}\,}

\def\id{\hbox{id}\,}
\def\Ad{\hbox{Ad}\,}

\def\Rep{{\bf Rep\,}}
\def\A{{\cal A}}
\def\F{{\cal F}}

\def\Hil{{\cal H}}

\def\D{{\cal D}}

\def\M{{\cal M}}

\def\one{{\bf 1}}
\def\onne{{\thinmuskip=5.5mu 1\!1\thinmuskip=3mu}}
\def\o{\otimes}

\def\cros{\raise1.9pt\hbox{$\scriptscriptstyle
          > $}\!\raise1.5pt\hbox{$\scriptstyle\triangleleft\,$}}

\def\C{\,{\raise1.5pt\hbox{$\scriptscriptstyle |$}
        \thinmuskip=4mu \!\!C\thinmuskip=3mu}}
\def\Z{{Z\!\!\!Z}}

\def\triv{D_{\varepsilon}}
\def\S{{\cal S}}
\def\Zro{{\cal Z}}
\def\amo{{\lower9pt\hbox{$\otimes$}\atop\raise2pt\hbox{
          $\scriptscriptstyle A^R\equiv\hat A^L$}}}

\def\w{\parbox{.22in}{\begin{picture}(10,10)(-5,-5)
       \put(0,0){\circle{6}}
       \end{picture}
       }}
\def\b{\parbox{.22in}{\begin{picture}(10,10)(-5,-5)
       \put(0,0){\circle*{6}}
       \end{picture}
       }}
\def\ww#1{\parbox{.7in}{\begin{picture}(42,23)(-5,-6)
          \put(3,3){\circle{6}}
          \put(6,3){\vector(1,0){17}}\put(21,3){\line(1,0){15}}
          \put(39,3){\circle{6}}
          \put(20,8){$#1$}
          \end{picture}
          }}
\def\wb#1{\parbox{.7in}{\begin{picture}(42,23)(-5,-6)
          \put(3,3){\circle{6}}
          \put(6,3){\vector(1,0){17}}\put(21,3){\line(1,0){15}}
          \put(39,3){\circle*{6}}
          \put(20,8){$#1$}
          \end{picture}
          }}

\def\bb#1{\parbox{.7in}{\begin{picture}(42,23)(-5,-6)
          \put(3,3){\circle*{6}}
          \put(6,3){\vector(1,0){17}}\put(21,3){\line(1,0){15}}
          \put(39,3){\circle*{6}}
          \put(20,8){$#1$}
          \end{picture}
          }}

\def\smallbbb#1#2#3#4{\parbox{.7in}{
    \begin{picture}(42,40)(-10,-15)
    \put(0,0){\circle*{6}}
    \put(3,0){\vector(1,0){14}}\put(17,0){\line(1,0){10}}
    \put(30,0){\circle*{6}}
    \put(14,-9){$#3$}
    \put(15,15){\circle*{6}}
    \put(2.25,2.25){\line(1,1){10}}\put(10,10){\vector(1,1){0}}
    \put(0,8){$#1$}
    \put(27.95,2.25){\line(-1,1){10}}\put(24,6){\vector(1,-1){0}}
    \put(24,8){$#2$}
    \put(12,4){$#4$}
    \end{picture}
    }}

\def\bigwbb#1#2#3#4{\parbox{0.9in}{
    \begin{picture}(42,40)(-10,-12)
    \put(0,0){\circle{6}}
    \put(3,0){\line(1,0){34}}\put(22,0){\vector(1,0){0}}
    \put(40,0){\circle*{6}}
    \put(18,-9){$#3$}
    \put(20,20){\circle*{6}}
    \put(2.25,2.25){\line(1,1){15}}\put(12,12){\vector(1,1){0}}
    \put(1,10){$#1$}
    \put(37.95,2.25){\line(-1,1){15}}\put(32,8){\vector(1,-1){0}}
    \put(34,10){$#2$}
    \put(17,6){$#4$}
    \end{picture}
    }}

\def\maunw#1#2#3#4#5{\parbox{.9in}{
    \begin{picture}(60,50)(-30,-25)
    \put(-20,0){\circle*{6}}
    \put(20,0){\circle*{6}}
    \put(0,20){\circle{6}}
    \put(0,-20){\circle{6}}
    \put(-17.5,2.5){\line(1,1){15}}
                     \put(-12,8){\vector(-1,-1){0}}
    \put(-17.5,-2.5){\line(1,-1){15}}
                     \put(-12,-8){\vector(-1,1){0}}
    \put(17.5,2.5){\line(-1,1){15}}
                     \put(12,8){\vector(1,-1){0}}
    \put(17.5,-2.5){\line(-1,-1){15}}
                     \put(12,-8){\vector(1,1){0}}
    \put(0,-17){\line(0,1){34}}\put(0,-8){\vector(0,-1){0}}
    \put(-20,9){$#1$}
    \put(-20,-17){$#2$}
    \put(15,9){$#3$}
    \put(15,-17){$#4$}
    \put(2,0){$#5$}
    \end{picture}
    }}
\def\maunb#1#2#3#4#5{\parbox{.9in}{
    \begin{picture}(60,50)(-30,-25)
    \put(-20,0){\circle{6}}
    \put(20,0){\circle{6}}
    \put(0,20){\circle*{6}}
    \put(0,-20){\circle*{6}}
    \put(-17.5,2.5){\line(1,1){15}}
                     \put(-7,13){\vector(1,1){0}}
    \put(-17.5,-2.5){\line(1,-1){15}}
                     \put(-7,-13){\vector(1,-1){0}}
    \put(17.5,2.5){\line(-1,1){15}}
                     \put(7,13){\vector(-1,1){0}}
    \put(17.5,-2.5){\line(-1,-1){15}}
                     \put(7,-13){\vector(-1,-1){0}}
    \put(0,-17){\line(0,1){34}}\put(0,-8){\vector(0,-1){0}}
    \put(-20,9){$#1$}
    \put(-20,-17){$#2$}
    \put(15,9){$#3$}
    \put(15,-17){$#4$}
    \put(2,0){$#5$}
    \end{picture}
    }}
\def\tetra#1#2#3#4#5#6{\parbox{.9in}{
    \begin{picture}(60,50)(-30,-25)
    \put(-20,0){\circle{6}}
    \put(20,0){\circle{6}}
    \put(0,20){\circle*{6}}
    \put(0,-20){\circle*{6}}
    \put(-17.5,2.5){\line(1,1){15}}
                     \put(-7,13){\vector(1,1){0}}
    \put(-17.5,-2.5){\line(1,-1){15}}
                     \put(-7,-13){\vector(1,-1){0}}
    \put(17.5,2.5){\line(-1,1){15}}
                     \put(7,13){\vector(-1,1){0}}
    \put(17.5,-2.5){\line(-1,-1){15}}
                     \put(7,-13){\vector(-1,-1){0}}
    \put(0,-17){\line(0,1){15}}\put(0,2){\line(0,1){14}}
                     \put(0,-10){\vector(0,-1){0}}
    \put(-17,0){\line(1,0){34}}
                     \put(10,0){\vector(1,0){0}}
    \put(-20,9){$#1$}
    \put(-20,-17){$#2$}
    \put(15,9){$#3$}
    \put(15,-17){$#4$}
    \put(-10,-2){$#5$}
    \put(0,6){$#6$}
    \end{picture}
    }}
\def\Wigner#1#2#3#4#5#6{\parbox{.9in}{
    \begin{picture}(60,50)(-30,-25)
    \put(-20,0){\circle{6}}
    \put(20,0){\circle*{6}}
    \put(0,20){\circle*{6}}
    \put(0,-20){\circle*{6}}
    \put(-17.5,2.5){\line(1,1){15}}
                     \put(-7,13){\vector(1,1){0}}
    \put(-17.5,-2.5){\line(1,-1){15}}
                     \put(-7,-13){\vector(1,-1){0}}
    \put(17.5,2.5){\line(-1,1){15}}
                     \put(12,8){\vector(1,-1){0}}
    \put(17.5,-2.5){\line(-1,-1){15}}
                     \put(7,-13){\vector(-1,-1){0}}
    \put(0,-17){\line(0,1){15}}\put(0,2){\line(0,1){14}}
                     \put(0,-10){\vector(0,-1){0}}
    \put(-17,0){\line(1,0){34}}
                     \put(10,0){\vector(1,0){0}}
    \put(-20,9){$#1$}
    \put(-20,-17){$#3$}
    \put(15,9){$#5$}
    \put(15,-17){$#6$}
    \put(-10,-2){$#2$}
    \put(0,6){$#4$}
    \end{picture}
    }}

\begin{document}
\large
\title{\bf A Coassociative C$^*$-Quantum Group\\[.15cm]
  \bf with Non-Integral Dimensions\\[1cm]}
\author{\sc Gabriella B\"ohm $^1$\\
\\
and\\
\\
\sc Korn\'el Szlach\'anyi $^2$\\
\\
Central Research Institute for Physics\\
H-1525 Budapest 114, P.O.B. 49, Hungary}
\date{September 7, 1995}

\maketitle

\footnotetext[1]{E-mail:
  BGABR@rmki.kfki.hu\\
  Supported by the Hungarian Scientific Research Fund, OTKA
  T 016 233}

\footnotetext[2]{E-mail:
  SZLACH@rmki.kfki.hu\\
  Supported by the Hungarian Scientific Research Fund,
  OTKA--1815.}

\vskip 2truecm

\begin{abstract}

By weakening the counit and
antipode axioms of a $C^*$-Hopf algebra and allowing for the
coassociative coproduct to be non-unital we obtain a
quantum group, that we call a {\em weak $C^*$-Hopf algebra},
which is sufficiently general to describe the
symmetries of essentially arbitrary fusion rules. This amounts to
generalizing the Baaj-Skandalis multiplicative unitaries to
multipicative partial isometries. Every weak $C^*$-Hopf algebra
has a dual which is again a weak $C^*$-Hopf algebra.
An explicit example
is presented with Lee-Yang fusion rules. We shortly discuss
applications to amalgamated crossed products, doubles, and
quantum chains.

\end{abstract}

\newpage
\normalsize

\sec{Introduction}

Conformal field theories provide examples of quantum field theory
models with finitely many superselection sectors $p,q,$\ldots
such that their intrinsic dimensions $d_p$ can take
non-integer values. Here $d_p$ denotes a component of the
Perron--Frobenius eigenvector of the fusion matrices,
$\sum_rN_{pq}^rd_r=d_pd_q$, but also coincides with the
statistical dimension in the sense of Algebraic QFT \cite{H}.
The symmetry of such a QFT cannot be described by a finite
dimensional $C^*$-Hopf algebra $H$. As a matter of fact if $H$ is
isomorphic to $\oplus_r \mbox{Mat}(n_r,\C)$ then $d_p=n_p$ solves
the above equation as a consequence of the unitalness of the
coproduct: $\Delta(\onne)=\onne\o\onne$. In order to fit to the
new situation, the concept of a Hopf algebra has been replaced by
more and more general structures:
quasi-Hopf algebras \cite{Dr}, truncated (or weak) quasi-Hopf
algebras \cite{MS2}, and rational Hopf algebras
\cite{V,FGV1}.

Although the rational Hopf algebra approach succesfully
reproduces all data of the quantum field theory encoded in the
representation category $\Rep\A$ of its observable algebra $\A$,
it has a serious flaw: Its dual is not an associative algebra,
therefore "crossed product" of $\A$ with a coaction of the
symmetry,
which is usually the algebra $\F$ of charge carrying fields, is
non-associative. The same non-associativity forbids also to define
an action of the symmetry on the fields without introducing at the
same time operators implementing the whole symmetry algebra
\cite{MS2}. In the theory of quantum chains one faces the
analogue problem if one would like to construct a quantum chain
with non-integer statistical dimensions. The strategy of the
construction \cite{NSz1,NSz2} requires the observable algebra to
be an iterated crossed
product $\dots\cros A\cros\hat A\cros A\cros\dots$ of a Hopf-like
algebra
$A$ and its dual $\hat A$. This would turn the observable algebra
to be non-associative as well as its "symmetry", the "double" of
$A$.

These problems naturally raise the question: Was it necessary
to give up coassociativity in order to
cover all interesting cases of
non-integral fusion rules? Our proposal of a {\it weak $C^*$-Hopf
algebra} shows that coassociativity can be maintained even in
the most general case.
Its axioms will be given in Sect.\,2. It
is a $C^*$-algebra $A$ and in applications to rational
quantum field theory models it is finite dimensional. In the
interesting cases it does not have, however, any 1-dimensional
representation. Therefore the counit $\varepsilon\colon A\to\C$
can
only be a coalgebraic counit and not an algebra map any more. The
coproduct $A\ni x\mapsto \Delta(x)\equiv x_{(1)}\o x_{(2)}$ is
coassociative but not unit preserving in
general. The antipode axiom is weakened as well so that the maps
$A\ni x\mapsto S(x_{(1)})x_{(2)}$ and
$A\ni x\mapsto x_{(1)}S(x_{(2)})$ become projections $\Pi^R$ and
$\Pi^L$, respectively, to certain non-trivial $C^*$-subalgebras of
$A$ that play important role in crossed products: These
subalgebras will give amalgamations (see Sect.\,3) between $A$ and
its dual $\hat
A$ in the Weyl algebras $A\cros\hat A$, $\hat A\cros A$, and in
the double $\D=A\bowtie\hat A$.

It is important to emphasize that our axioms define a selfdual
algebraic structure, just like the axioms for a Hopf algebra. That
is the dual $\hat A$ can be given structural maps
$\hat\Delta,\hat\varepsilon,\hat S$ satisfying the same axioms as
$\Delta,\varepsilon$, and $S$ do for $A$.

The simplest example for such a weak Hopf algebra can be found
[see Appendix] in
studying the quasi-double $\D^{\omega}(G)$ of a finite group $G$
\cite{DiPaRo}. $\D^{\omega}(G)$ has a 1-dimensional
block corresponding to the trivial representation but its
coproduct is quasi coassociative. "Blowing up" the double by a
full matrix algebra $M_n$, that is introducing the algebra
$\M^{\omega}(G):=\D^{\omega}(G)\o M_n$,
where $n$ can be taken to be the
order of $G$, there is a coassociative coproduct on
$\M^{\omega}(G)$ producing precisely the fusion rules $N_{pq}^r$
of $\D^{\omega}(G)$. This coassociative coproduct is related to
the quasi-coassociative coproduct of $\D^{\omega}(G)$ by a
"skrooching" (a name proposed by J. Stasheff \cite{St}), i.e. by
conjugation with a
partial isometry $U\in\M^{\omega}(G)\o\M^{\omega}(G)$.
Transforming out a non-trivial cocycle $\omega$ by skrooching
was not possible in $\D^{\omega}(G)$ but there is more flexibility
in the blown up version $\M^{\omega}(G)$.
Not only their fusion rules coincide but also
the representation category of $\M^{\omega}(G)$ is
equivalent to that of the quasi double. It turns out that
$\M^{\omega}(G)$ itself is a kind of double. It is the double of a
blown up version of the Hopf algebra of complex functions on $G$.

The dimensions of $\M^{\omega}(G)$ are of course integers.
In order to find examples with non-integer dimensions one needs
to construct weak $C^*$-Hopf algebras $A$ with $S^2\neq\id$.
Namely, one can prove that the dimensions $d_p$ of a weak
$C^*$-Hopf algebra are all integers if the antipode is involutive.
(Recall that in case of compact matrix pseudogroups \cite{W}
$S^2=\id$
follows if $A$ is finite dimensional. Therefore non-integer
dimensions can be expected only for infinite dimensional compact
matrix pseudogroups.) There is a general construction of finite
dimensional weak $C^*$-Hopf algebras from a solution $F$ of the
pentagon equation, typically having non-integral $d_p$-s and
therefore $S^2\neq\id$.
One starts from a given set of fusion matrices $N_{pq}^r$ $^3$
\footnotetext[3]{Here we restrict ourselves to multiplicity free
fusion rules: $N_{pq}^r\leq 1$.}
and a corresponding
solution $F(^{pqr}_{\ s})_{tu}$ of the pentagon equations
\beq
\sum_e \ F(^{pqr}_{\ s})_{ef}\ F(^{aer}_{\ d})_{cs}\
 F(^{apq}_{\ c})_{be}\ =\ F(^{apf}_{\ d})_{bs}\ F(^{bqr}_{\
d})_{cf}
\eeq
and the unitarity conditions
\beq
\sum_e\ \overline{F(^{pqr}_{\ s})}_{ef}\ F(^{pqr}_{\
s})_{eg}\ = \ \delta_{fg}\ N_{qr}^gN_{pg}^s\ .
\eeq
Such a solution is provided in any quantum field theory with
finitely many superselection sectors as it was emphasized in
\cite{MS1,S}. In spite of its tempting similarity to the
pentagon
equation for the reassociator $\varphi$ of a quasi-Hopf algebra
\cite{Dr} triviality of $\varphi$ is not the same as triviality of
$F$.
Hence a coassociative quantum symmetry (trivial $\varphi$) does
not imply any kind of simplification on the solution $F$. Quite on
the contrary, we will see in Sect.4 that eqn. (1.1) can always be
read as the coassociativity condition of an appropriate coproduct.

The weak $C^*$-Hopf algebra $A$ we construct from (1.1--2)
has fusion rules precisely
given by $N_{pq}^r$ and its $6j$-symbols by $F(^{pqr}_{\
s})_{tu}$. As a $C^*$-algebra it is isomorphic
to $\oplus_pM_{n_p}$ where $n_p=\sum_{ab}N_{ap}^b$. For example in
case of the Lee-Yang fusion rules one obtains (see Sect.\,5)
the algebra $A=M_2\oplus
M_3$ which can be thought of an inhomogeneous "blowing up" of the
quasi-coassociative rational Hopf algebra $M_1\oplus M_2$ of
\cite{FGV1}. For the Ising fusion rules our
construction yields $A=M_3\oplus M_3\oplus M_4$ instead of the
quasi-coassociative $M_1\oplus M_1 \oplus M_2$.

We briefly mention that weak $C^*$-Hopf algebras can also be
described
by the methods of Baaj and Skandalis \cite{BaSk}. Let $\Hil$ be a
Hilbert space and $V$ a partial isometry acting on $\Hil\o\Hil$.
We call $V$ a {\it multiplicative partial isometry} if
the following two equations hold true on $\Hil\o\Hil\o\Hil$.
\beanon
V_{12}V_{13}V_{23}&=&V_{23}V_{12}\\
V_{13}V_{23}V_{23}^*&=&V_{12}^*V_{12}V_{13}
\eeanon
To this data there is an associated pair $(A,\hat A)$ of weak
$C^*$-Hopf algebras in duality. The reverse statement can be
verified easily. Namely, if $A$ is a weak $C^*$-Hopf algebra then
$A$ is dense in the Hilbert space $\Hil$ carrying the
left regular representation of $A$ and
\[V(x\o y)\ :=\ x_{(1)}\o x_{(2)}y\qquad x,y\in A\subset\Hil\]
defines a multiplicative
partial isometry. Its initial and final
projections are $V^*V=\Delta(\onne)$ and
$VV^*=\hat\Delta(\hat\onne)$, respectively.

Most of our results refer, at present, to finite dimensional weak
$^*$-Hopf algebras although there seems to be no
obstruction towards a generalization
to the infinite dimensional case.

Ideas about a coassociative $C^*$-quantum group very similar to
ours have been proposed earlier by F. Nill
(unpublished).

\sec{Axioms of Weak $C^*$-Hopf Algebras}

\begin{defi}
A weak $^*$-Hopf algebra is a $^*$-algebra
$A$ with unit $\onne$ together with linear maps $\Delta\colon A\to
A\o A$, $\varepsilon\colon A\to \C$, and $S\colon A\to A$ called
the coproduct, the counit, and the antipode respectively, if
the following axioms hold:
\begin{axioms}
\Delta(xy)&=&\Delta(x)\Delta(y)\\
\Delta(x^*)&=&\Delta(x)^*\\
(\Delta\o\id)\circ\Delta&=&
              (\id\o\Delta)\circ\Delta
\end{axioms}
\begin{axioms}
\varepsilon(xy)&=&\varepsilon(x\onne_{(1)})
                    \varepsilon(\onne_{(2)}y)\\
\varepsilon(x^*x) &\geq& 0\\
(\varepsilon\o\id)\circ\Delta\ =\ &\id&
           =\ (\id\o\varepsilon)\circ\Delta
\end{axioms}
\begin{axioms}
S(xy) &=& S(y)S(x)\\
S\circ^*\circ\ S\ \circ^* &=& \id\\
\Delta\circ S &=& (S\o S)\circ \Delta^{op}
\end{axioms}
\begin{axiom}
S(x_{(1)})x_{(2)}\o x_{(3)} = \onne_{(1)}\o
x\onne_{(2)}
\end{axiom}
for all $x, y\in A$. If furthermore $A$ possesses a faithful
$^*$-representation on a Hilbert space it is called a weak
$C^*$-Hopf algebra.
\end{defi}

If we add to these axioms either the condition
$\Delta(\onne)=\onne\o\onne$, or the condition
$\varepsilon(xy)=\varepsilon(x)\varepsilon(y)$, or the condition
$S(x_{(1)})x_{(2)}=\onne\varepsilon(x)$
then $A$ becomes a usual
$(C)^*$-Hopf algebra. We do not know whether
a simpler set of axioms exists that are equivalent to (A.1--4).
But our axioms go definitely beyond Hopf algebras
by the vast class of examples with non-integer dimensions
constructed in Sect.\,4.

Now we give, without proofs, some of the consequences of the above
axioms. Details will be published elsewhere [BNSz]. The maps
$\Pi^{L/R}$ defined by
\beq
\Pi^L(x):=x_{(1)}S(x_{(2)})\,,\qquad\Pi^R(x):=S(x_{(1)})x_{(2)}
\eeq
are not conditional expectations, are not even
$^*$-preserving, but are linear
projections that project onto $C^*$-subalgebras $A^L$,
resp. $A^R$ of $A$. They are isomorphic as $C^*$-algebras, the
antipode maps $A^L$ onto $A^R$, $A^L$ commutes with $A^R$,
and the restriction of $S^2$ onto $A^{L/R}$ is the identity.
The restriction of the counit onto $A^L$ (or $A^R$) is a faithful
trace. We may then define an orthonormal basis $\{e^i\}$ of $A^L$,
i.e. $\varepsilon(e^{i*}e^j)=\delta^{i,j}$. It follows that
$e_i:=S(e^{i*})$ defines an orthonormal basis for $A^R$. The
coproduct of the unit can be expressed as $\Delta(\onne)=e_i\o
e^i$ (summation over $i$ is understood). It follows that
$N:=\varepsilon(\onne)$ is a positive integer equal to the
dimension of $A^{L/R}$.

The dual $\hat A$ of $A$ is defined to be the space of linear
functionals $\varphi$ on $A$ and is equipped with a multiplication
and
a comultiplication obtained by tranposing the comultiplication and
multiplication of $A$ w.r.t the canonical pairing $\langle\ ,\
\rangle\colon \hat A \times A\to\C$. The unit element of $\hat A$
is $\hat\onne:=\varepsilon$. The antipode
$\hat S$ and the star operation of $\hat A$ are defined
respectively by
\bea
\langle\hat S(\varphi),x\rangle &=& \langle\varphi,S(x)\rangle\\
\langle\varphi^*,x\rangle &=& \overline{\langle
\varphi,S(x)^*\rangle}
\eea
\begin{thm} {\sl If $(A,\onne,\Delta,\varepsilon,S)$
satisfies
Axioms (A1--4) then $(\hat A,\hat\onne,\hat\Delta,\hat\varepsilon,
\hat S)$ satisfies Axioms (A1--4), too. That is the notion of a
weak $^*$-Hopf algebra is selfdual.}
\end{thm}
The proof requires a little work because our axioms are not so
transparently selfdual as those of a Hopf algebra. For example
in order to verify the crucial axioms (A.2a) and (A.4) in the dual
one proves at first that the relations
\beanon
\onne_{(1)}\o\varepsilon(\onne_{(2)}\onne_{(1)'})\onne_{(2)'}&=&
\onne_{(1)}\o\onne_{(2)}\\
S(x_{(1)})x_{(2)}y&=&y_{(1)}\varepsilon(xy_{(2)})\quad x,y\in A
\eeanon
are consequences of (A1--4). Then pairing the first with
$\varphi\o\psi\in \hat A\o \hat A$ and the second with
$\varphi\in\hat
A$ one obtains the dual counterparts of (A.2a) and (A.4),
respectively.

Just like in case of a pair of dual Hopf algebras there are
natural left and right actions of $A$ on $\hat A$
given by the Sweedler's arrows
\beq
x\rightarrow\varphi:=\varphi_{(1)}\langle\varphi_{(2)},x\rangle\,,
\qquad \varphi\leftarrow
x:=\langle\varphi_{(1)},x\rangle\varphi_{(2)}
\eeq
for all $\varphi\in\hat A$ and $x\in A$. Similarly one defines the
left and right actions of $\hat A$ on $A$.

The canonical pairing $\langle\ ,\ \rangle$ restricted to $\hat
A^L\times A^R$ is non-degenerate and so is its restriction to
$\hat A^R\times A^L$. The bases dual to $\{e_i\}$ and $\{e^i\}$
are found to be $\hat e^i=\hat\onne\leftarrow e^i\in\hat A^L$ and
$\hat e_i=e_i\to\hat\onne\in\hat A^R$, respectively.

{}From now on we assume that $A$ satisfies (A.1--4) and is a
finite dimensional $C^*$-algebra ("finite quantum group").
Let $\Rep A$ denote the category of finite dimensional (not
necessarily non-degenerate) $^*$-representations of $A$. Its
arrows are the intertwiners $T\colon D_1\to D_2$ with the
unit arrow at the object $D$ being $\one_D=D(\onne)$.
$\Rep A$ becomes a monoidal category if we define for objects
$D_1,\
D_2$ the product representation by $D_1\times D_2:=(D_1\o D_2)
\circ\Delta$ and for arrows $T\colon D_1\to D_2$, $T'\colon
D'_1\to D'_2$ the product arrow by $T\times T':=(T\o T')\cdot
(D_1\times D'_1)(\onne)$. The unit object (i.e. the trivial
representation of $A$) is defined to be the GNS representation
$\triv$ associated to the state ${1\over N}\cdot\varepsilon$.
Using the identity
\beq
\varepsilon(x^*y)\ =\ \varepsilon(\Pi^L(x)^*\Pi^L(y))
\eeq
one recognizes that the representation space of $\triv$ is $A^L$
equipped with the Hilbert-Schmidt scalar product
$(x^L,y^L):=\varepsilon(x^{L*}y^L)$, $x^L,y^L\in A^L$. The action
is given by $\triv(x)\Pi^L(y)=\Pi^L(xy)$, $x,y\in A$. Sometimes it
is convenient to use matrix elements of $\triv$ in the orthonormal
basis $\{e^{i*}\}$ of $A^L$: $\triv^{ij}(x)=\varepsilon(e^ixe_j)$.
The name trivial representation for $\triv$ is justified by the
existence of isometric intertwiners $u^L_D\colon D\to\triv\times
D$ and $u^R_D\colon D\to D\times\triv$ that
satisfy the triangle identities
\bealph
u_{D_1\times D_2}^L &=& u_{D_1}^L\times\one_{D_2}\\
\one_{D_1}\times u_{D_2}^L &=& u_{D_1}^R\times\one_{D_2}\\
u_{D_1\times D_2}^R &=& \one_{D_1}\times u_{D_2}^R
\eealph
and are natural in $D$. They can be defined by the matrix
elements
\bea
(u_D^L)^{i\alpha,\beta} &:=& D^{\alpha\beta}(e^i)\\
(u_D^R)^{\alpha i,\beta} &:=& D^{\alpha\beta}(e_i^*)
\eea

The usual formula $\bar D=D^T\circ S$ for the conjugate
representation does not work unless $S^2=\id$, since $\bar D$ is
not unitary in general. However, there exists an $S_0\in\End A$
and an invertible $C\in A$ such that $S_0^2=\id$, $S_0\circ\ ^*=\
^*\circ S_0$, and $S=\Ad_C\circ S_0$. Then we may define the
conjugate of the representation $D$ by
\beq
\bar D^{\alpha\beta}(x)\ :=\ D^{\beta\alpha}(S_0(x))\qquad x\in
A\,.
\eeq
Rigidity intertwiners can be obtained as
follows. Let
\beq
(R_D)^{\bar\alpha\alpha,i}\ :=\
D^{\alpha\bar\alpha}(e_iS_0(C))
\eeq
then $R_D$ is an injective intertwiner from
$\triv$ to $\bar D\times D$. Furthermore let $\bar R_D:=R_{\bar
D}$. Then the rigidity relations
\bea
u_D^{L*}(\bar R_D^*\times\one_D)(\one_D\times R_D)u_D^R &=&
\one_D\\
u_{\bar D}^{R*}(\one_{\bar D}\times\bar R_D^*)
(R_D\times\one_{\bar D})u_{\bar D}^L &=& \one_{\bar D}
\eea
hold provided we have adjusted $C$ to satisfy $S(C^*)C=\onne$,
which is always possible in view of the freedom in multiplying $C$
by a central invertible element. We have therefore

\begin{thm} If $A$ is a finite dimensional weak $C^*$-Hopf algebra
then the category $\Rep A$ of its finite dimensional
$^*$-representations is a rigid monoidal $C^*$-category.
\end{thm}
\noindent Coassociativity of the coproduct leads to "almost"
strict monoidality of $\Rep A$. Although the reassociator maps
$\varphi_{D_1,D_2,D_3}\colon(D_1\times D_2)\times D_3\to
D_1\times(D_2\times D_3)$ are the trivial identity arrows one
needs the non-trivial isometric arrows $u_D^L$ and $u_D^R$ to
compare $D$ with $\triv\times D$ or $D\times\triv$.

Somewhat surprisingly the trivial representation need not be
irreducible. Consider the example where $A$ is a finite
dimensional Abelian $C^*$-algebra with 1-dimensional irreducible
representations $D_p$ obeying the fusion rules $D_p\times
D_q=\delta_{p,q}D_q$. Now it is clear that the "trivial"
representation is not only reducible but even faithful:
$\triv=\sum_pD_p$. However, this quantum group is a sum of
(1-dimensional) quantum groups each of them having an irreducible
trivial representation. This is a general phenomenon as the next
Proposition shows.
\begin{prop} For $A$ a finite dimensional weak
$C^*$-Hopf algebra let $\S$ denote the set of its irreducible
sectors and let $\triv=\bigoplus_{p\in\S}\
\oplus_{\mu=1}^{\nu_p}\
D_p$ be the decomposition of the trivial representation into
irreducibles.

Then $\nu_p\leq 1$, i.e. $\triv$ is multiplicity free,
$\bar p=p$ for all $p\in\S$ such that $\nu_p=1$, i.e. the trivial
representation contains only selfconjugate irreducibles.
Let $\Zro:=\{p\in\S|\nu_p=1\}$ and for $p\in\Zro$ define the
subsets of sectors
 \[  \S_p:=\{q\in\S|N_{pq}^q=1\}\ .   \]
Then $\{\S_p|p\in\Zro\}$ is a partition of $\S$ each member of
which is closed under the monoidal product (i.e. $p\in\Zro,\
q_1,q_2\in\S_p$, and $N_{q_1q_2}^r>0$ imply that $r\in\S_p$) and
under conjugation ($p\in\Zro,\ q\in\S_p$ imply $\bar q\in\S_p$).
The irreducible $D_p$ serves as a monoidal unit for the full
subcategory $\Rep_p A$ of $\Rep A$ generated by the irreducibles
$q\in\S_p$. If $p\neq p'$, the monoidal product of any $D\in
\Rep_p A$ and $D'\in \Rep_{p'} A$ is the zero object, $D\times D'=0$,
and the intertwiner space $(D|D')$ is zero dimensional.
\end{prop}

\begin{coro} Any finite dimensional weak $C^*$-Hopf
algebra $A$ can be decomposed uniquely into a sum $A=\oplus_p A_p$
of weak $C^*$-Hopf subalgebras such that the counits
$\varepsilon_p=\varepsilon|_{A_p}$ of $A_p$ are pure
(unnormalized) states.
\end{coro}

If $A$ has a counit that is pure $A$ will be called {\em pure}.
For pure $A$ the rigidity intertwiners can be normalized to obey
\beq
R_D^*R_D=d_D\one_{\triv}\,,\qquad \bar R_D^*\bar
R_D=d_D\one_{\triv}  \eeq
with a uniquely determined positive number $d_D$, called
the dimension of $D$. The function $D\mapsto d_D$ is additive and
multiplicative for direct sums and for the monoidal product of
representations.

Like compact groups finite dimensional weak $C^*$-Hopf algebras
possess
unique Haar measures in the following sense. There exists a unique
$h\in A$ characterized by the property
\beq
xh=\Pi^L(x)h\,,\quad hx=h\Pi^R(x)\,,\qquad\forall x\in A
\eeq
and by the normalization conditions $\Pi^L(h)=\onne=\Pi^R(h)$.
For Hopf algebras $\Pi^{L/R}=\varepsilon$, therefore this
definition coincides with the usual one \cite{Sw}.
The Haar measure $h$ satisfies the following important properties:
\bealph
\langle\varphi^*\varphi,h\rangle &\geq& 0\qquad \varphi\in\hat
A\\
\varphi\geq 0,\ \langle\varphi,h\rangle=0 &\Rightarrow&
\varphi=0\\
\langle\hat\onne,h\rangle &=& N\\
h^2\ =\ h&=&h^*\ =\ S(h)\\
 h_{(1)}x\o h_{(2)} &=& h_{(1)}\o h_{(2)}S(x)\qquad x\in A
\eealph
It follows that the Haar integral $\int\varphi dh\equiv
\langle\varphi,h\rangle$ is left and right invariant:
\bea
\langle (x\rightarrow\varphi)\psi,h\rangle &=&
\langle \varphi(S(x)\rightarrow\psi) ,h\rangle\\
\langle (\varphi\leftarrow x)\psi,h\rangle &=&
\langle \varphi(\psi\leftarrow S^{-1}(x)) ,h\rangle
\eea
Furthermore the maps
\beq
E^L(\varphi):=h\to \varphi\,,\quad E^R(\varphi):=
\varphi\leftarrow h      \label{cond exp}
\eeq
define conditional expectations from $\hat A$ onto $\hat A^L$ and
$\hat A^R$,
respectively. In order to show its existence we give an explicit
expression for $h$.
Choose the orthonormal basis $\{e^i\}$ to be diagonal with
respect to the decomposition $A=\oplus A_p$ of Corollary 2.5. Then
$\triv^{ij}\equiv \hat e^i\hat e_j$ is zero unless $i$ and $j$
belong to the same block $p\in\Zro$. Choose matrix units
$\{e_p^{ij}\}$ in the trivial block of $A_p$ that are dual to
$\triv^{ij}$, i.e. $\langle\triv^{ij},e_p^{kl}\rangle=\delta^{ik}
\delta^{jl}\delta_{k\in p}\delta_{l\in p}$ . Then
\beq
h\ =\ \sum_{p\in\Zro}{1\over N_p}\sum_{i,j\in p}
\varepsilon(e^i)\,e_p^{ij}\,\varepsilon(e_j)
\eeq
and its coproduct takes the following form using matrix units
$\{e_r^{\alpha\beta}\}$, $\alpha,\beta=1,\dots,n_r$ of $A$:
\beanon
\Delta(h) &=& (S\o\id)(X)\cdot(C\o
C^{-1})\\
 \hbox{where}\quad X &=& \sum_{p\in\Zro}{1\over N_p}
\sum_{r\in\S_p}{1\over d_r}\sum_{\alpha,\beta=1}^{n_r}
\ e_r^{\alpha\beta}\o e_r^{\beta\alpha}\ =\ X^*\ =\ X^{op}\ .
\eeanon
Since the existence of the Haar measure in $A$
implies that $\hat A\ni\varphi\mapsto\langle\varphi,h\rangle$
is a faithful positive linear functional, $\hat A$ is also a
$C^*$-algebra. This shows that selfduality in the sense of Theorem
2.2 is also true for finite dimensional $C^*$-weak Hopf algebras.

Like in \cite{W} one can prove that the Hermitean invertible
element $s:=CC^*$ implementing the square of the antipode,
$sxs^{-1}=S^2(x),\ x\in A$, is grouplike, that is
$\Delta(s)=(s\o s)\Delta(\onne)=\Delta(\onne)(s\o s)$.
Although our finite quantum group is always unimodular in the
sense that it has a 2-sided Haar measure, the modular operator in
the left regular representation is non-trivial because the Haar
state is not a trace:
\beq
\langle\varphi\psi,h\rangle=
\langle\psi(s\rightarrow\varphi\leftarrow s),h\rangle
\eeq
If $A$ is pure the dimension function can be expressed as
$d_D={1\over N}\,\chi_D(s^{-1})$, where $\chi_D$ denotes ordinary
(unnormalized) trace in the representation $D$.

\sec{The Weyl algebras and the Double}

The Weyl algebra $A\cros\hat A$ of a weak $^*$-Hopf algebra $A$ is
defined to be the unital $^*$-algebra generated by $A$ and $\hat
A$, as unital $^*$-subalgebras, subjected to the commutation
relation
\beq
\varphi x\ =\ x_{(1)}\,\langle
x_{(2)},\varphi_{(1)}\rangle\,\varphi_{(2)}
\eeq
It follows that the Weyl algebra is a $^*$-algebra in
which the Sweedler arrows $\hat A\to A$ and $\hat A\leftarrow A$
become inner:
\bea
\varphi\rightarrow x &=& \varphi_{(1)}\,x\,\hat S(\varphi_{(2)})\\
\varphi\leftarrow x &=& S(x_{(1)})\,\varphi\,x_{(2)}
\eea
If $A$ is a Hopf $^*$-algebra then its Weyl algebra is isomorphic
to $A\o \hat A$ as a linear space and is simple as an algebra [N].
In our more general setting $A\cros\hat A$ is an amalgamated
tensor product $A\amo\hat A$ over the common $^*$-subalgebra
$A^R\cong\hat A^L$ identified through the isomorphism $A^R\ni
x^R\mapsto (\hat\onne\leftarrow x^R)\in\hat A^L$. In other words
$\hat e^i=e_i^*$ is an identity in $A\cros \hat A$ for all
$i=1,\dots,N$. The other two "edge" subalgebras turn out to be the
relative commutants of $A$ and $\hat A$ within the Weyl algebra,
thus we have
\beanon
A\cap\hat A\ =\ A^R&\equiv&\hat A^L\\
A'\cap(A\cros\hat A) &=& \hat A^R\\
(A\cros\hat A)\cap \hat A ' &=& A^L
\eeanon
One can define the other Weyl algebra $\hat A\cros A$ by
interchanging in (3.1) the roles of $\varphi$ and $x$.

If $A$ is a weak $C^*$-Hopf algebra then using the Haar measure
$h$ one can make $\hat A$ to be a Hilbert space with the scalar
product $(\varphi,\psi):=\langle\varphi^*\psi,h\rangle$. The left
regular representation of $\hat A$ thus becomes a faithful
$^*$-representation $\pi$, which is nothing but the GNS
representation associated to the Haar state
$\langle\cdot,h\rangle$. Now left and right invariance of the Haar
measure implies that $\pi$ can be extended both to a
$^*$-representation $\pi^R$ of $A\cros\hat A$ and to a
$^*$-representation $\pi^L$ of $\hat A\cros A$:
\[
\parbox{2in}{
\beanon
\pi^L(x)\psi&=&x\to\psi\\
\pi^L(\varphi)\psi&=&\varphi\psi
\eeanon }
\parbox{2in}{
\beanon
\quad\pi^R(x)\psi&=&\psi\leftarrow S^{-1}(x)\\
\quad\pi^R(\varphi)\psi&=&\varphi\psi
\eeanon}
\]
Both representations turn out to be faithful. Hence Weyl
algebras of weak $C^*$-Hopf algebras are $C^*$-algebras.

The projections $\pi^L(h)$ and $\pi^R(h)$ in the left regular
representation project onto the subspaces $\hat A^L$ and $\hat
A^R$, respectively. Now using standard technics
\cite{J,GHJ} it is
not difficult to show that the two Weyl algebras arise through the
basic construction from the inclusions $\hat A^R\subset\hat A$ and
$\hat A^L\subset \hat A$, respectively, in such a way that
$\pi^L(h)$ and $\pi^R(h)$ become the Jones projections with the
associated conditional expectations being given by (\ref{cond
exp}). Utilizing
also the dual statements we obtain the following Jones triples
\[\begin{array}{lr}
A^{L}\subset A\subset A\cros\hat A\quad&
\quad\hat A\cros A\supset A\supset A^{R}\\
\hat A^{L}\subset\hat A\subset \hat A\cros A\quad&
\quad A\cros\hat A\supset\hat A\supset\hat A^{R}
\end{array}\]

The double $\D=A\bowtie\hat A$ of a weak $^*$-Hopf algebra $A$ can
be defined
by generators and relations exactly in the same way as for Hopf
algebras (e.g. see Appendix B of \cite{NSz1}). $\D$ is a unital
$^*$-algebra generated
by symbols $\D(x),\ x\in A$ and $\D(\varphi),\ \varphi\in\hat A$.
The relations require that $\D(A)$ and $\D(\hat A)$ form unital
subalgebras of $\D$ isomorphic to $A$ and $\hat A$, respectively.
Furthermore the following commutation relation is postulated:
\beq
\D(\varphi)\D(x)=\D(x_{(2)})\D(\varphi_{(2)})\,\langle\varphi_{(1)},
x_{(3)}\rangle\,\langle\varphi_{(3)},S^{-1}(x_{(1)}\rangle
\eeq
for all $x\in A,\ \varphi\in\hat A$. It follows that $\D$ is an
amalgamated tensor product of $A$ and $\hat A$ with the common
$^*$-subalgebras being $\partial A:=<A^L,A^R>$ of $A$ and
$\partial\hat A:=<\hat A^L,\hat A^R>$ of $\hat A$ with
identifications given between $A^R$ and $\hat A^L$ as in the Weyl
algebra $A\cros\hat A$ and between $\hat A^R$ and $A^L$ as in the
Weyl algebra $\hat A\cros A$. In other words the following
identities hold in $\D$:
\beq
\D(e^i)=\D(\hat e_i^*),\quad \D(e_i)=\D(\hat e^{i*}),\qquad
i=1,\dots,N\,.
\eeq

One can introduce a coproduct $\Delta_{\D}$, a counit
$\varepsilon_{\D}$, and an antipode $S_{\D}$ on the double just
like in case of Hopf algebras.
These structural maps turn out to satisfy all axioms (A.1--4).
Hence the double of a weak $^*$-Hopf algebra is a weak $^*$-Hopf
algebra again.

Similar statement holds also in the $C^*$ case. As a matter
of fact let $A$ be a weak $C^*$-Hopf algebra and $h$, $\hat h$
denote the
Haar measures in $A$ and $\hat A$, respectively. Then the map
$\hat\D\ni\Phi\mapsto
\langle \Phi,\D(h)\D(\hat h)\rangle$
can be shown to define a faithful positive linear functional.
Thus $\hat\D$, and therefore its dual $\D$ too, is
a weak $C^*$-Hopf algebra.

Like for Hopf algebras the double is always quasitriangular. That
is if $\{b^k\}$ and $\{\beta_k\}$ denote bases of $A$ and $\hat A$
that are in duality, $\langle\beta_k,b^l\rangle=\delta_k^l$, then
\beq{\cal R}\ :=\ \sum_k\ \D(b^k)\o\D(\beta_k)\quad \in\
\D\o\D\eeq
is a partial isometry satisfying the usual hexagon identities and
intertwines between $\Delta_{\D}$ and the opposite coproduct
$\Delta_{\D}^{op}$.

A quantum chain can now be constructed using methods of
\cite{NSz1,NSz2}
based on a finite dimensional weak $C^*$-Hopf algebra $A$. Its
observable algebra is defined to be the infinite crossed product
$\A=\dots A\cros\hat A\cros A\cros\hat A\cros\dots$
equipped with the following
local structure. For each interval
$I=\{i,i+1,\dots,j\}\subset\Z$ we have the local algebra
$\A(I)\equiv\A_{i,j}:=
A_i\cros A_{i+1}\cros\dots\cros A_j$ with the 1-point algebras
being $A_{2i}=A$, $A_{2i+1}=\hat A$. The 2-point algebras
$\A_{i,i+1}$ in this net are the Weyl algebras. It follows that
the net will fail to be split, i.e. $\A(I)$ will not be simple for
any $I$, and will not have the intersection property, i.e.
$A_i\cap A_{i+1}\neq\C\onne$ (provided $A^L$ is not simple). The
failure of these two properties can encourage one to hope that
the statistical dimensions will be non-integers. Indeed, if we
define a localized coaction $\rho\colon \A\to\A\o\D$ of the double
$\D$ of $A$ by formula (3.18) of \cite{NSz1} then we exhibit the
double as a quantum symmetry of the superselection sectors
(whether it is also universal is not clear at present). Therefore
the statistical dimensions and other representation theoretic data
can be determined solely from the knowledge of the double. The
example treated in Sect. 5 will demonstrate that the double can
really have non-integral dimensions.

Crucial is for the construction that the net $\A$ satisfies Haag
duality. Like for Hopf algebras Haag duality in $\A$ is intimately
related to the existence of integrals in $A$ and $\hat A$.
As a matter of fact if $h_i$ denotes the unique Haar integral in
$A_i$ then the maps
\beq
\eta^L_i(a):=h_{i(1)}\,a\,S(h_{i(2)})\ ,\qquad
\eta^R_i(a):=S(h_{i(1)})\,a\,h_{i(2)}\ ,\qquad a\in\A
\eeq
define conditional expectations satisfying a generalized
$\eta$-property (cf. \cite{NSz1}):
\beq
\eta^L_i(\A_{-\infty,i-1})=
\A_{-\infty,i-2}\,,\qquad
\eta^R_i(\A_{i+1,\infty})=
\A_{i+2,\infty}
\eeq
which implies, as in \cite{NSz1}, Haag duality of the
net $\A$.

Field algebras can be constructed by taking crossed products with
respect to actions $\rho$ of the dual of the double:
$\F=\A\cros_{\rho}\hat\D$ in which a local subalgebra of
$\A$ becomes amalgamated with the "left subalgebra" $\hat\D^L$.
Although the inclusion $\A\subset\F$ in general is not
irreducible, $\A'\cap\F=\hat\D^R$, the gauge principle is
satisfied. Namely,
$\F$ carries a natural action $\gamma$ of the symmetry
algebra $\D$ such that $\A$ reappears in $\F$ as the
invariant subalgebra: $\gamma_{h_{\D}}(\F)=\A$, where
$h_{\D}$ is the Haar integral of $\D$. In other words, in the
double crossed product $\A\cros_{\rho}\hat\D\cros_{\gamma}\D$ the
commutant of $\D$ is the observable algebra $\A$.

\sec{Weak $C^*$-Hopf Algebras and the Pentagon Equation}

As it was pointed out in \cite{MS1} any quantum field theory
determines
a solution $F$ of the pentagon equations (1.1) which in turn
coincides with the recoupling coefficients (6j-symbols) of the
underlying quantum symmetry. The task is then to reconstruct a
quantum group from its recoupling coefficients $F$. Here we shall
sketch an argument how a weak $C^*$-Hopf algebra can be
reconstructed from a solution $F$ of the pentagon equations.
(The braiding structure which is also an important ingredient in
QFT is neglected throughout this paper.)

We formulate the argument in a quite general combinatorial way
which was strongly motivated by Ocneanu's quantum cohomology
\cite{O}.
Let $K$ be a finite simplicial 3-complex having only 2 vertices:
\w and \b. There are 3 types of (oriented) edges:
\bb p, \wb i, and \ww{\hat p}, but no edges go from
\b to \w. Faces may be attached to three edges $a,b,c$
if their endpoints fit into a figure of the kind
\bigwbb abc{\nu} although possibly with some other (allowed)
coloring of the vertices.
The faces bordered by the triangle $(a,b,c)$ are labelled by
$\nu=1,\dots,N^c_{ab}$. (If all the three edges $p,q,r$ connect
\b with \b then the number of such faces $N^r_{pq}$ can
be thought the fusion coefficients of a quantum field theory or of
a monoidal category, although we do not want to make any
restriction
on these numbers.) One 3-simplex is attached to four
faces whenever they form the boundary of a tetrahedron.

Now we construct two $C^*$-algebras: $A$ associated to the vertex
\b and $\hat A$ associated to \w. Each edge $\bb q$
corresponds to a simple block of $A$.
The simple block $q$ is by definition $\End V_q$ where $V_q$ is
the Hilbert space an orthonormal basis of which is given by the
triangles \bigwbb iqj{\nu} with all possible values of $i,j,$ and
$\nu$. In order to facilitate the notation from now on we assume
that $N_{ab}^c\leq 1$. Then a matrix unit basis for $A$ is
provided by the elements
\beq
e_q^{(i'j')(ij)}\ =\ \left(^{i'\,i}_{j'j}\right)_q\ \equiv\
\maunb{i'}{j'}ijq\label{e}
\eeq
Obviously $A$ is a matrix algebra $\oplus_q M_{n_q}$ with
dimensions $n_q=\sum_{i,j}N_{iq}^j$.

One defines $\hat A$ in the same way but the roles of \w and
\b being interchanged. Hence $\hat A\cong\oplus_{\hat q}
M_{m_{\hat
q}}$ where for an edge $\ww {\hat q}$ the block dimension $m_{\hat
q}$ is equal to
$\sum_{i,j} N_{\hat qj}^i$. Matrix units of $\hat A$ are therefore
\beq
\hat e_{\hat q}^{(i'j')(ij)}\ =\
\left[^{i'\,i}_{j'j}\right]_{\hat q}
\ \equiv\ \maunw{i'}{j'}ij{\hat q}\label{ehat}
\eeq

The coproduct on $A$ is defined by means of a $\C$-valued 3-chain
$F_1$
\[
F_1(^{iq'q''}_{\ k})_{jq}\ =\ \Wigner ijkq{q'}{q''}
\]
supported on tetrahedrons with one \w and three
\b as their vertices.
\[
\Delta\left(\left(^{i'\,i}_{j'j}\right)_q\right)\ :=\
\sum_{q'q''kk'}\ \left(^{i'\,i}_{k'k}\right)_{q'}\otimes
\left(^{k'k}_{j'j}\right)_{q''}\ \cdot\
F_1(^{i'q'q''}_{\ j'})_{k'q}\ \overline{F_1(^{iq'q''}_{\ j})}
_{kq}
\]
This coproduct will be a coassociative $^*$-algebra map provided
$F_1$ satisfies the unitarity condition (1.2) and furthermore
there exists a function $F_0$ supported on the tetrahedra with
no vertices of type \w such that the pentagon equation
\[
\sum_e \ F_0(^{pqr}_{\ s})_{ef}\ F_1(^{aer}_{\ d})_{cs}\
 F_1(^{apq}_{\ c})_{be}\ =\ F_1(^{apf}_{\ d})_{bs}\ F_1(^{bqr}_{\
d})_{cf}
\qquad\qquad(P_1)\]
holds and $F_0$ also satisfies unitarity. One recognizes in $F_1$
the analogue of
the Wigner-coefficients or $3j$-symbols while $F_0$ comprises the
Racah-coefficients. Thus $(P_1)$ is nothing but the recoupling
equation for the quantum group $A$. The notation $(P_1)$ refers
to that among the 5 vertices present in this equation there is
exactly 1 of type \w. Of course for the existence of an
$F_0$ satisfying $(P_1)$ it is necessary that $F_0$ satisfies a
pentagon equation of the type $F_0F_0F_0=F_0F_0$, which is
precisely eqn.(1.1), also called the Elliot--Biederharn relation.
In our notation it is $(P_0)$ because it
involves no \w. Similarly, in order to have a coproduct on
$\hat A$ one postulates the existence of a 3-chain $F_3$ (the
"dual Wigner-coefficients") supported on tetrahedra with 3
\w and a 3-chain $F_4$ (the "dual Racah-coeffients")
supported on tetrahedra with no \b at all. They
satisfy the dual recoupling equation $(P_4)$, the dual
pentagon equation $(P_5)$, and unitarity conditions.

In order to construct a (non-degenerate bilinear) pairing between
$A$ and $\hat A$ we use the suggestion coming from geometry. A
matrix unit (\ref{e}) of $A$ is a pair of triangles with a common
edge
of type \bb~ with the other two vertices being white. A matrix
unit of $\hat A$ on the other hand is a pair of triangles
(\ref{ehat}) with
\ww~ as the common edge and the other two vertices being black.
Hence rotating one of them by $180^{\circ}$ around its "main
diagonal" we can glue them
together to obtain a tetrahedron with exactly 2 white vertices.
Hence the pairing will be defined in terms of a 3-chain $F_2$
(an "Ocneanu cell") so that
\beq
\langle\left[^{i'\,i}_{j'j}\right]_{\hat p}\ ,
\ \left(^{k'k}_{l'\,l}\right)_q \rangle\ =\
\delta_{i'k'}\delta_{j'k}\delta_{i,l'}\delta_{j,l}
\cdot\tetra{k'}{l'}kl{\hat p}q
\eeq
The condition for this pairing to transpose the coproduct of $A$
to the product of $\hat A$ is again a pentagon of the type
$(P_2)$: $F_2F_1F_1=F_2F_2$. Similarly, the compatibility of the
product on
$A$ with the coproduct on $\hat A$ is the pentagon equation
$(P_3)$: $F_2F_3F_3=F_2F_2$.

Summarizing, if $F=(F_0,F_1,F_2,F_3,F_4)$ is a unitary 3-cocycle
on $K$, i.e. satisfies the Big Pentagon equation
$(P)=(P_0,P_1,P_2,P_3,P_4,P_5)$ and unitarity, then its
restrictions $F_n$ $n=0,\dots,4$
to tetrahedra with fixed coloring of the vertices determine the
various structural maps of a weak $C^*$-Hopf algebra $A$. We
recall that once having the $C^*$-algebras $A$ and $\hat A$ in
duality the antipode can be introduced by the formula
$S:=\ ^*\circ\ _*$ where the antilinear involution $_*$ is defined
by
$\langle\varphi,x_*\rangle:=\overline{\langle\varphi^*,x\rangle}$.
The resulting weak $C^*$-Hopf algebra will be pure iff the complex
$K$ was 2-connected. For such pure quantum groups $A$ we have
derived from the axioms the existence of a unique block $q=0$, the
trivial representation, a unique involution $q\mapsto\bar q$
of the blocks of $A$ determined by the action of the antipode
on the center, and one such involution $\hat
p\mapsto\bar{\hat p}$ of the dual blocks. Therefore the fact that
the fusion ring determined by the fusion coefficients $N_{ab}^c$
is associative, unital, and involutive are not assumptions on the
complex $K$ but all are consequences of the existence of a
solution for the Big Pentagon equation.

In practice, however, we have only given a solution $F_0$ of
$(P_0)$. That is $F$ is known only on the subcomplex $K_A\subset
K$ containing the single vertex \b. There seems to be a canonical
way to extend $F_0$ to a solution $F$ determining a selfdual weak
$C^*$-Hopf algebra $A$, i.e. such that $\hat A\cong A$ as weak
$C^*$-Hopf algebras. In order to have such an extension, however,
we need to postulate --- on $K_A$ at least --- the existence of
the trivial block 0 and the involution $q\mapsto \bar q$ with
the usual properties
\bealph
&&N_{0b}^c=\delta_{bc}=N_{b0}^c\\
&&N_{ab}^0=\delta_{a\bar b}\\
&&N_{c\bar b}^a=N_{ab}^c=N_{\bar a c}^b
\eealph
Now the extension goes as follows. Define $K$ by adjoining to
$K_A$ another copy of it denoted $K_{\hat A}$ with its vertex
denoted by \w and then connect the two vertices by as many
edges as $K_A$ has. Therefore all three types of edges of $K$ will
carry the same label $p$ running over the set of sectors (of $A$
or of the underlying quantum field theory). Triangles with mixed
vertices are set between three edges $(p,q,r)$ whenever
there is a triangle between the edges of the same labels in $K_A$.
With this choice
we obtain $A\cong\hat A\cong\oplus_p M_{n_p}$ where
$n_p=\sum_{q,r}
N_{pq}^r$. Now the extension $F$ is obtained by setting
$F_1=F_3=F_4=F_0$ and $F_2$
to be the inverse of $F_0$ in the sense of the equation
\beq
\sum_r\ F_2(^{pqr}_{\ s})_{tu}\ F_0(^{p'qr}_{\ s})_{tu}\
=\ \delta^{pp'}\cdot\hbox{constraints}
\eeq
for each fixed value of $q,s,t,u$. The
existence of this kind of inverse (and therefore of the pairing)
follows since
(4.4abc) imply that $F_0$ has an $S_4$ symmetry \cite{FGV2} hence
unitarity implies its invertibility also in the required indices.

Therefore a selfdual unimodular finite quantum group can be
associated to any rational quantum field theory describing its
internal symmetry on the level of representation categories.
Of course, this does not solve the problem of how to determine
unique integers $n_q$
\cite{S}, i.e. the problem of uniqueness of the quantum symmetry.
There may be
other extensions $K$ of $K_A$ and solutions $F$ on $K$ of the Big
Pentagon equation extending the given $F_0$. Thus there may be
other weak $C^*$-Hopf algebras $B$ (though perhaps not selfdual)
such that $\Rep B$ is equivalent to $\Rep A$.

\sec{An Example with Lee--Yang Fusion Rules}

Consider the fusion ring generated by two selfconjugate sectors
$0$ and $1$ where $0$ is the
unit element and $1\times 1=0+1$. Corresponding to these fusion
rules there is a solution $F$ of the pentagon equation
\cite{FGV1}.

The weak $C^*$-Hopf algebra $A_{L.Y.}$ given below has been
constructed from that solution using the general method of
Sect.\,4.
Setting the free parameter $\zeta$ in \cite{FGV1} to 1 all
structural
maps can be written as integer polynomials in $z=\sqrt
{(\sqrt5-1)/2}$. The intrinsic dimensions of the two sectors are
$d_0=1$ and $d_1=z^{-2}\equiv(\sqrt 5+1)/2$.

As a $C^*$-algebra $A_{L.Y.}=M_2\oplus M_3$. We fix matrix units
$e_0^{ij}$ in $M_2$ and $e_1^{ij}$ in $M_3$. The coproduct
is given by
\beanon
\Delta(e_0^{11})&=&e_0^{11}\o e_0^{11}+e_1^{11}\o
e_1^{33}\\
\Delta(e_0^{12})&=&e_0^{12}\o e_0^{12}+z^2e_1^{13}\o
e_1^{31}+ze_1^{12}\o e_1^{32}\\
\Delta(e_0^{22})&=&e_0^{22}\o e_0^{22}+z^4e_1^{33}\o e_1^{11}+
z^3e_1^{32}\o e_1^{12}+z^3e_1^{23}\o e_1^{21}+z^2e_1^{22}\o
e_1^{22}\\
\Delta(e_1^{11})&=&e_0^{11}\o e_1^{11}+e_1^{11}\o
e_0^{22}+e_1^{11}\o e_1^{22}\\
\Delta(e_1^{12})&=&e_0^{12}\o e_1^{12}+e_1^{12}\o e_0^{22}+
ze_1^{13}\o e_1^{21}-z^2e_1^{12}\o e_1^{22}\\
\Delta(e_1^{13})&=&e_0^{12}\o e_1^{13}+e_1^{13}\o e_0^{21}
+e_1^{12}\o e_1^{23}\\
\Delta(e_1^{22})&=&e_0^{22}\o e_1^{22}+e_1^{22}\o e_0^{22}\\
&&+z^2e_1^{33}\o e_1^{11}-z^3e_1^{32}\o e_1^{12}-z^3e_1^{23}\o
e_1^{21}+z^4e_1^{22}\o e_1^{22}\\
\Delta(e_1^{23})&=&e_0^{22}\o e_1^{23}+e_1^{23}\o e_0^{21}
+ze_1^{32}\o e_1^{13}-z^2e_1^{22}\o e_1^{23}\\
\Delta(e_1^{33})&=&e_0^{22}\o e_1^{33}+e_1^{33}\o
e_0^{11}+e_1^{22}\o e_1^{33}\\
\eeanon
The counit and the antipode are as follows:
\beanon
\varepsilon(e_0^{ij})&=&1\quad i,j=1,2,\qquad
\varepsilon(e_1^{ij})=0\quad i,j=1,2,3\,,\\
S(e_0^{ij})&=&e_0^{ji}\quad i,j=1,2\qquad
S(e_1^{ij})=z^{i-j}\,e_1^{\bar j\,\bar i}\quad i,j=1,2,3\,,
\eeanon
where we introduced the notation $\bar 1=3,\ \bar 2=2,\ \bar
3=1$. One can check easily that they satisfy the axioms (A1--4).
The left and right subalgebras turn out to be Abelian, isomorphic
to $M_1\oplus M_1$, with orthonormal bases given by the minimal
projections
\beanon
&&e^1=e_0^{11}+e_1^{11},\qquad\qquad
e^2=e_0^{22}+e_1^{22}+e_1^{33}\\
&&e_1=e_0^{11}+e_1^{33},\qquad\qquad
e_2=e_0^{22}+e_1^{11}+e_1^{22}
\eeanon
The projections $\Pi^{L/R}$ take the following simple form
\beanon
  \Pi^L(e_0^{ij})=e^i&\qquad &\Pi^R(e_0^{ij})=e_j\\
  \Pi^L(e_1^{ij})=0&\qquad &\Pi^R(e_1^{ij})\,=\,0
\eeanon
The Haar measure $h$ is a rank 1 projection belonging to the
support of $\varepsilon$, i.e. of the trivial block 0:
\[h\ =\ {1\over 2}\sum_{i,j=1,2}\ e_0^{ij}\,.\]

The above weak $C^*$-Hopf algebra is obviously the smallest in
dimension among the weak $C^*$-Hopf algebras with $S^2\neq \id$.
Nevertheless computing the structural maps of its double
is quite a horrible task. Since the existence of a double with
non-integer intrinsic dimensions is crucial for the construction
of quantum chains with non-integer statistical dimensions, we have
calculated --- using some computer aid --- the block structure,
the fusion rules, and the dimensions $d_p$ of the double
$\D_{L.Y.}$
of the Lee-Yang quantum group $A_{L.Y.}$. The results are the
following. As an algebra $\D_{L.Y.}$ is isomorphic to
\[\D_{L.Y.}\ =\ M_2\oplus M_3\oplus M_3\oplus M_5\]
All the four sectors are selfconjugate.
Denoting these sectors by $2,3,3',5$, respectively we have the
following fusion rules. 2 is the trivial representation, therefore
$2\times p=p\times 2=p$ for all sectors $p$. The fusion is
commutative: $p\times q=q\times p$ for all $p,q$. Furthermore
\[\begin{array}{lcr}
 3\times 3\,=2+3\ \ &\,\ 3'\times 3'=2+3'\ &\ \ 5\times
                                                5=2+3+3'+5\\
 3\times 3'=5\ \ &\ 3'\times 5\,=3+5&\\
 3\times 5\,=3'+5
 \end{array}\]
These fusion rules yield the intrinsic dimensions $d_2=1,\
d_3=d_{3'}=z^{-2},\
d_5=z^{-4}$. Except $d_2$ these are not integers and are related
to the $d_1$ of $A_{L.Y.}$ in a very
simple way, $d_2=d_1^0,\ d_3=d_1,\ d_{3'}=d_1,\ d_5=d_1^2$ which
may have its explanation in the selfduality of $A_{L.Y.}$.

\vskip 2truecm
{\bf\noindent Acknowledgement}: One of us (K. Sz.) wishes to thank
Florian Nill for his suggestions and remarks on the notions of
integral and $C^*$-structure.

\vskip 2truecm
{\twelvebf\noindent Appendix: Blowing up the Quasi-Double
$\D^{\omega}(G)$}
\bigskip

Let $G$ be a finite group of order $|G|=N$ and $\omega$ a
$U(1)$ valued 3-cocycle
on $G$. Then $\D^{\omega}(G)$ is the quasi-Hopf algebra with
product and coproduct given on the basis elements $(g,h)\in
G\times G$ as follows \cite{DiPaRo}:
\beanon
(g,h)(g',h')&=&\delta_{gh,hg'}\,\theta_g(h,h')\cdot(g,hh')\\
\mbox{where}&&\theta_g(h,h')={
\omega(h,h',h'^{-1}h^{-1}ghh')\,\omega(g,h,h')
\over\omega(h,h^{-1}gh,h')}\\
\Delta_{\omega}((g,h))&=&\sum_{k\in G}\gamma_h(k,k^{-1}g)\cdot
(k,h)\o (k^{-1}g,h)\\
\mbox{where}&&\gamma_h(x,y)={
\omega(h,h^{-1}xh,h^{-1}yh)\,\omega(x,y,h)
\over\omega(x,h,h^{-1}yh)}
\eeanon
$\Delta_{\omega}$ is quasi-coassociative and no skrooched version
$\Delta'_{\omega}=\Ad_u\circ\Delta_{\omega}$ of it can be
coassociative, unless $\omega$ is a coboundary.

Consider the full matrix algebra $M_N$ with a fixed choice of
matrix units $\{e_{ab}\}$ and make it a weak $C^*$-Hopf algebra
by introducing $\Delta_M(e_{ab})=e_{ab}\o e_{ab}$,
$\varepsilon_M(e_{ab})=1$, $S_M(e_{ab})=e_{ba}$ for
$a,b=1,\dots,N$. Now we claim that on the tensor product
$\M^{\omega}(G)=\D^{\omega}(G)\o M_N$, which is a weak quasi-Hopf
algebra, there is a skrooching
$U\in\M^{\omega}(G)\o\M^{\omega}(G)$ transforming the
quasi-coassociative non-unital coproduct
$\Delta_{\omega}\o\Delta_M$ into a coassociative non-unital
coproduct $\Delta$:
\beanon
\Delta(\xi\o m)&=&U\cdot(\Delta_{\omega}\o\Delta_M)
(\xi\o m)\cdot U^*\\
U&=&\sum_{a,b,c\in G}\omega^{-1}(a,b,c)\,[(c,1)\o
e_{ab,a}]\o [(b,1)\o e_{a,a}]\ .
\eeanon
The counit remains the tensor product
$\varepsilon=\varepsilon_{\omega}\o\varepsilon_M$ while the
antipode undergoes a skrooching $S=\Ad_V\circ(S_{\omega}\o S_M)$
where $V=U_1\cdot (S_{\omega}\o S_M)(U_2)$.
With the structure maps $(\Delta,\varepsilon,S)$ the blown up
double $\M^{\omega}(G)$ becomes a weak $C^*$-Hopf algebra.
It is also quasitriangular with R-matrix $U^{op}(R_1\o I\o R_2 \o
I)U^*$ where $R=R_1\o R_2$ denotes the old R-matrix of the
quasi-double.
This blown up double and the original quasi-Hopf algebra
have identical representation theories, namely
$\Rep\M^{\omega}(G)$ and
$\Rep\D^{\omega}(G)$ are equivalent as braided monoidal
categories. $\M^{\omega}(G)$ depends non-trivially on
the cohomology class of $\omega$ although it is coassociative for
any choice of the 3-cocycle.

$\M^{\omega}(G)$ can be shown to be the double in the sense of
Sect.\,3 of a weak $C^*$-Hopf algebra $A^{\omega}(G)$. The latter
one can be defined as a blowing up of the Hopf algebra $\C(G)$ of
complex functions on $G$. If $\delta_g,\ g\in G$ denote the
minimal
projections in $\C(G)$ then $A^{\omega}(G)$ is obtained from the
tensor product weak Hopf algebra $\C(G)\o M_N$ by skrooching with
the partial isometry
\[ u=\sum_{a,b,c\in G}\omega^{-1}(a,b,c)\,\left[\delta_b\o
e_{a,a}\right]\o\left[\delta_c\o e_{ab,a}\right]\ .
\]
Equivalently $A^{\omega}(G)$ can be constructed by the method of
Sect.\,4 from the cocycle $\omega$ itself, if it is interpreted as
a solution $F=\omega^{-1}$
of the pentagon equation on an appropriate complex $K_A=K(G)$.
As a matter of fact, the complex $K(G)$ has one vertex \b and $N$
edges labelled by $g\in G$. A face is attached to the triangle
\smallbbb ghk~ if and only if the zero curvature condition $gh=k$
is satisfied. The resulting $C^*$-algebra $A=A^{\omega}(G)$ is
isomorphic to $\C(G)\o M_N$ and its Wigner and Racah coefficients,
as well as its Ocneanu cell, are given by the cocycle $\omega$.
The relation between the two constructions can be elucidated by
the formula
\[\delta_g\o e_{ab}\ =\ \maunb a{ag}b{bg}g\]

\eject

\end{document}